\input{style/arxiv-general.cfg}
\documentclass[MSNbibl,number,citesort,seceqn,dvips]{arxbj}
\makeatletter
   \@ifpackageloaded{graphicx}{}{\usepackage{graphicx}}
\makeatother
\usepackage{upgreek}
\usepackage{graphicx}

\aid{0}
\volume{21}
\issue{2}
\pubyear{2015}
\firstpage{1047}
\lastpage{1066}
\doi{10.3150/14-BEJ597} 

\makeatletter

\newcommand{\rright}{\right}
\newcommand{\lleft}{\left}
\newcommand{\rrvert}{\vert}
\newcommand{\llvert}{\vert}
\newtheorem{theorem}{Theorem}[section]
\newproclaim{assumption}{Assumption}[section]
\newtheorem{lemma}{Lemma}[section]
\newtheorem{corollary}{Corollary}[section]
\newremark{example}{Example}[section]
\newtheorem{proposition}{Proposition}[section]

\newcommand{\eE}{\mathrm{E}}
\newcommand{\eqref}[1]{(\ref{#1})}
\newcommand{\Cov}{\operatorname{Cov}}
\newcommand{\diag}{\operatorname{diag}}
\newcommand{\Var}{\operatorname{Var}}
\newcommand{\Corr}{\operatorname{Corr}}
\renewcommand{\pi}{\uppi}

\makeatother

\begin{document}
\begin{frontmatter}

\title{Functional partial canonical correlation}
\runtitle{Functional partial canonical correlation}

\begin{aug}
\author[A]{\inits{Q.}\fnms{Qing} \snm{Huang}\corref{}\thanksref{e1}\ead[label=e1,mark]{qing.huang@asu.edu}} \and
\author[A]{\inits{R.}\fnms{Rosemary} \snm{Renaut}\thanksref{e3}\ead[label=e3,mark]{Rosie.Renaut@asu.edu}}
\address[A]{School of Mathematics and Statistics, Arizonas State
University, Tempe, AZ 85287, USA.\\
\printead{e1,e3}}
\end{aug}

\received{\smonth{5} \syear{2012}}
\revised{\smonth{12} \syear{2013}}

%
\begin{abstract}
A rigorous derivation is provided for canonical correlations and
partial canonical correlations for certain Hilbert space indexed
stochastic processes. The formulation relies on a key congruence
mapping between the space spanned by a second order,
$\mathcal{H}$-valued, process and a particular Hilbert function space
deriving from the process' covariance operator. The main results are
obtained via an application of methodology for constructing orthogonal
direct sums from algebraic direct sums of closed subspaces.
\end{abstract}

%
\begin{keyword}
\kwd{congruent Hilbert space}
\kwd{covariance operator}
\kwd{Hilbert space indexed process}
\kwd{orthogonal direct sum}
\end{keyword}

\end{frontmatter}

\section{Introduction}
\label{sec1}
Canonical correlation analysis (CCA) is one of the principal tools for
studying the relationship between two random vectors in multivariate
analysis. There have now been several attempts to widen the definition
of CCA to include vectors of infinite length and, more generally,
stochastic processes (see, e.g., Eubank and Hsing \cite{eandh} and
references
therein). Functional canonical correlation falls into this latter
category wherein one obtains data that represent the sample paths of
continuous time processes. In this paper we provide a framework for
canonical correlation and partial canonical correlation analysis for a
class of stochastic processes that includes those arising in
functional data.

A somewhat general formulation assumes that we have a probability
space $(\Omega,\mathcal{A},P)$, a real, separable Hilbert space
$\mathcal{H}$, with norm and inner product $\Vert\cdot\Vert$ and
$\langle\cdot,\cdot\rangle$ and an $\mathcal{H}$-valued random variable $X$ in the
sense of Laha and Rohatgi \cite{landr}; that is, $X\dvtx \Omega
\rightarrow\mathcal{H}$ is a
measurable function relative to the Borel $\sigma$-field generated by
the class of all open subsets of $\mathcal{H}$. Our attention will be
restricted to random variables with $\mathrm{E}\Vert X\Vert^2 <\infty
$ with
expectation being relative to $P$. Associated with such a random
variable we can define the Hilbert space indexed
process
%
\begin{equation}
\label{eqhsi} Z(f) = \langle X,f\rangle
\end{equation}
for $f\in\mathcal{H}$. Then, from Vakhania \textit{et al.} \cite{vak} there
exists a mean
element $h \in\mathcal{H}$
and a covariance operator $S$ such that $\mathrm{E}[\langle X,f\rangle ]=\langle h,f\rangle $ and
$\mathrm{E}[\langle X-h,f\rangle \langle X-h,f'\rangle ]=\langle f,Sf'\rangle $ for all $f$, $f' \in
\mathcal{H}$. For simplicity, we assume that $\Vert h\Vert=0$. In
that case,
the covariance operator is determined by
%
\begin{equation}
\label{ecovop} \mathrm{E}\bigl[\langle X,f\rangle \bigl\langle X,f'
\bigr\rangle \bigr]=\bigl\langle f, Sf'\bigr\rangle .
\end{equation}

It is well known that $S$ in \eqref{ecovop} is a trace class
operator and therefore admits the eigenvalue--eigenvector decomposition
%
\begin{equation}
\label{eexpansion} S = \sum_{j = 1}^{\infty}
\lambda_j\phi_j\otimes\phi_j,
\end{equation}
where $\lambda_1 \geq\lambda_2 \geq\cdots>0$ are the eigenvalues,
$\phi_j$ is the eigenvector associated with $\lambda_ j$ and
$(f\otimes g)h = \langle f, h \rangle g$ for $f, g, h \in
\mathcal{H}$. A suitably normed version of the range of $S$ gives us
the reproducing kernel Hilbert space
%
\begin{eqnarray} \label{eHS}
\mathcal{H}(S)= \Biggl\{f\dvt f=\sum_{j=1}^\infty
\lambda_{j}f_{j}\phi_{j}, \Vert f
\Vert^2_{\mathcal{H}(S)}=\sum_{j=1}^\infty
\lambda _{j}f_{j}^2 = \bigl\Vert S^{-1/2}f
\bigr\Vert^2 < \infty \Biggr\}
\end{eqnarray}
that includes $\mathcal{H}$ as a proper subset when $S$ is not finite
dimensional which we hereafter assume to be the case.
The reproducing kernel Hilbert space recasts the range of $S$ under a
weaker norm where $S$ is invertible, since the Picard condition
(Engl \textit{et al.} \cite{engl})
\begin{eqnarray*}
\sum_{j = 1}^{\infty}\frac{\langle f, \phi_j\rangle^2 }{\lambda
_j} = \sum
_{j = 1}^{\infty}\lambda_j
f^2_j < \infty
\end{eqnarray*}
is satisfied for $f \in\mathcal{H}(S)$.
For each $f \in
\mathcal{H}(S)$ there corresponds a random variable
\begin{eqnarray*}
Z(f) = \sum_{j = 1}^{\infty}f_j
\langle X, \phi_j \rangle.
\end{eqnarray*}
These types of random variables are well defined and include those in
the process (\ref{eqhsi}) as a special case. Thus, for inferential
purposes we can focus on the Hilbert space
%
\begin{equation}
L_{Z}^2= \Biggl\{Z(f) = \sum
_{j=1}^\infty f_{j}\langle X,
\phi_j\rangle \dvt \bigl\Vert Z(f)\bigr\Vert_{L_Z^2}^2 :=
\Var\bigl(Z(f)\bigr)=\sum_{j=1}^{\infty}
\lambda_{j}f_{j}^2 < \infty \Biggr\}
\label{eLZ}
\end{equation}
which consists of all the linear combinations of the $\langle X,\phi_j\rangle $ that
have finite variance. Note that in addition to serving as an index
set, $\mathcal{H}(S)$ is isometrically isomorphic or congruent to
$L^2_Z$: a relationship that will be exploited in the
sequel. Parzen \cite{parzen3} calls $\mathcal{H}(S)$ a
congruent reproducing
kernel Hilbert space.

For functional data, $X$ and the $\phi_j$ are typically functions on
some continuous index set $E$. In that instance it follows from
Kupresanin \textit{et al.} \cite{kup} that working with $L^2_Z$ is equivalent
to working with the
space spanned by the $X$ process: that is,
%
\begin{equation}
\label{eql2subx} L_{X}^2 = \overline{ \Biggl\{a\dvt a = \sum
_{j=1}^na_{j}
X(t_j), t_j \in E, a_{j} \in\mathbb{R}, n =
1, 2, \ldots \Biggr\}}
\end{equation}
under the inner product $\eE[ab]$ for $a, b \in L_{X}^2$. In
fact, functional canonical correlation can be treated directly from
this latter perspective using reproducing kernel Hilbert space
techniques along the lines of those employed in Eubank and Hsing
\cite{eandh}. However,
our present formulation in terms of $L^2_Z$ has certain advantages (both
mathematical and computational) and appears to generalize more readily
to deal with partial canonical correlation and related ideas.

Assume now that we have two $\mathcal{H}$-valued random variables
$X_i, i = 1, 2$, whose associated covariance operators $S_i, i = 1,
2$, have the eigenvalue--eigenvector sequences $\{(\lambda_{ij},
\phi_{ij})\}_{j = 1}^{\infty}$ from \eqref{eexpansion}. These, in
turn, produce Hilbert spaces $L_{Z_i}^2, i = 1, 2$, defined analogous
to \eqref{eLZ} for processes $Z_i(f_i), i = 1,2$, that are indexed by
Hilbert spaces $\mathcal{H}(S_i)$ defined as in \eqref{eHS}. Then,
the (first) canonical correlation between $Z_1$ and $Z_2$ is defined
to be
%
\begin{equation}
\label{eccadef} \rho^2 = \sup_{\Vert f_i\Vert_{\mathcal{H}(S_i)} = 1, i = 1, 2}
\Cov^2\bigl(Z_1(f_1), Z_2(f_2)
\bigr).
\end{equation}
One can deduce from Eubank and Hsing \cite{eandh} that \eqref
{eccadef} is well defined
with the supremum being attained. We provide an independent
verification of this fact in the next section. If $\tilde{f}_1, \tilde
{f}_2$ are
maximizing functions, then $Z_1(\tilde{f}_1), Z_2(\tilde{f}_2)$ are
the first
canonical variables of the $Z_1$ and $Z_2$ processes,
respectively. Subsequent canonical correlations and variables can be
obtained similar to the first in an iterative process that parallels
the one employed in the standard multivariate analysis case; see,
for example, Eubank and Hsing \cite{eandh}.

A number of articles dealing with functional canonical correlation and
related concepts have focused on the case where the $Z_i(f_i)$ are
restricted to have
%
\begin{equation}
\label{eqsilly} \sum_{j = 1}^{\infty}f_{ij}^2
< \infty,\qquad  i = 1, 2,
\end{equation}
which has the consequence that $\sum_{j=1}^{\infty}f_{ij}\phi_{ij}
\in
\mathcal{H}$. In such instances the supremum \eqref{eccadef} need not
be attained as demonstrated in Cupidon \textit{et al.} \cite{cupid1} and
Cupidon \textit{et al.} \cite{cupid2}.
Dauxois and Pousse \cite{f1}, Dauxois \textit{et al.} \cite{f2}, Dauxois and Nkiet \cite{f3} and Dauxois \textit{et al.} \cite{f4} largely ignore this
issue with the consequence that their statistical applications become
relevant only for finite dimensional covariance operators whose ranges
are necessarily closed. Such results are, of course, already subsumed
by the original Hotelling \cite{hotel} work. In contrast,
He \textit{et al.} \cite{he} impose
restrictions on the cross-covariances of coefficients in the two
processes' Karhunen--Lo\`{e}ve expansions to insure that
(\ref{eqsilly}) is satisfied. Such restrictions are unnecessary as
will be seen in the next section.

In the present paper, we are interested not only in functional CCA but
functional partial canonical correlation, as well. In the case of
finite dimensional covariance operators, the idea was proposed by
Roy \cite{roy}. Given three random vectors $X_1, X_2$
and $X_3$, the
partial canonical correlation of $X_2$ and $X_3$ relative to $X_1$ was
defined as the ordinary canonical correlation between $\tilde{X}_2 =
X_2 - P_{X_1}X_2$ and $\tilde{X}_3 = X_3 - P_{X_1}X_3$, where
$P_{X_1}$ denotes projection onto the linear space spanned by
$X_1$. Related work by Dauxois and Nkiet \cite{f3} and Dauxois \textit{et al.}
\cite{f4} comes with the
restriction of a closed range for covariance operators which, again,
confines statistical applications to the finite dimensional setting
that was already treated in Roy's original work. In Section~\ref{sec3}, we show
how the partial canonical correlation concept can be rigorously
extended to infinite dimensions and functional data.

In the next section, we set out the main ideas that are needed for
rigorous treatment of canonical correlation and related concepts in
the context of Hilbert space indexed processes of the basic form
\eqref{eLZ}. The driving force behind our approach is the isometry
that exists between the $L_{Z}^2$ and $\mathcal{H}(S)$ spaces. To
demonstrate the utility of this analytic framework, we illustrate the
idea with two processes in the next section and extend this to three
processes and partial canonical correlation in Section~\ref{sec3}.

\section{CCA}
\label{sec2}
In this section, we begin with the case of two processes and establish
the properties of canonical correlations and variables as defined in
\eqref{eccadef}. Most of the basic techniques that are needed for the
three process setting of the next section are illustrated in this
somewhat simpler scenario thereby making it the natural starting point
for our exposition.

As in Section~\ref{sec1}, assume that we have two
$\mathcal{H}$-valued random variables with associated covariance
operators $S_i, i = 1, 2$, having eigenvalue--eigenvector
sequences $\{(\lambda_{ij},\phi_{ij})\}_{j = 1}^{\infty}$. From
Vakhania \textit{et al.} \cite{vak}, it may be concluded that there are also
cross-covariance
operators $S_{12}$ and $S_{21}$ defined by, for example,
\begin{eqnarray*}
\mathrm{E}\langle X_1,f_1\rangle \langle
X_2,f_2\rangle =\langle f_1,S_{12}f_2
\rangle
\end{eqnarray*}
with $S_{21}=S_{12}^*$ for $S_{12}^*$ the adjoint of $S_{12}$.

Now we construct a new Hilbert space
\begin{eqnarray*}
\mathcal{H}_0= \Biggl\{h=(f_1,f_2)\dvt
f_i\in\mathcal{H}(S_i), i=1,2, \Vert h
\Vert_0^2=\sum_{i = 1}^2
\Vert f_i\Vert^2_{\mathcal{H}(S_i)} <\infty \Biggr\}
\end{eqnarray*}
from which we obtain the $\mathcal{H}_0$ indexed process
\begin{eqnarray*}
Z(h)=Z_1(f_1)+Z_2(f_2)
\end{eqnarray*}
with covariance function
%
\begin{eqnarray}\label{ecovZh}
\Cov\bigl(Z(h),Z\bigl(h'\bigr)\bigr)&=& \Cov\bigl(Z_1(f_1),Z_1
\bigl(f_1'\bigr)\bigr)+\Cov\bigl(Z_2(f_2),Z_2
\bigl(f_2'\bigr)\bigr)
\nonumber
\\
&&{}+\Cov\bigl(Z_1(f_1),Z_2
\bigl(f_2'\bigr)\bigr)+\Cov\bigl(Z_1
\bigl(f_1'\bigr),Z_2(f_2)\bigr)
\nonumber
\\[-8pt]\\[-8pt]
&=& \bigl\langle f_1,f_1'\bigr\rangle
_{\mathcal{H}(S_1)}+\bigl\langle f_2,f_2'\bigr
\rangle _{\mathcal{H}(S_2)}
\nonumber
\\
&&+\Cov\bigl(Z_1(f_1),Z_2
\bigl(f_2'\bigr)\bigr)+\Cov\bigl(Z_1
\bigl(f_1'\bigr),Z_2(f_2)\bigr).\nonumber
\end{eqnarray}
In order to avoid the degenerate setting where perfect prediction is
possible, we impose the following condition.

%
\begin{assumption}
\label{aa1}
There exist no $(f_1,f_2) \in\mathcal{H}_0$ such that
$|\Corr(Z_1(f_1),Z_2(f_2))|=1$.
\end{assumption}

The cross-covariance terms in \eqref{ecovZh} can be characterized as
deriving from operators between $\mathcal{H}(S_1)$ and
$\mathcal{H}(S_2)$. To see this, define the functional
\begin{eqnarray*}
l_{f_2}(f_1)=\Cov\bigl(Z_1(f_1),Z_2(f_2)
\bigr)
\end{eqnarray*}
on $\mathcal{H}(S_1)$. Clearly, $l_{f_2}$ is linear since covariance
is bilinear and, e.g., $Z_1(\alpha f_1+\alpha'f_1') = \alpha
Z_1(f_1)+\alpha'Z_1(f_1')$ for any scalars $\alpha$, $\alpha'$ and any
$f_1$, $f_1' \in\mathcal{H}(S_1)$. Also, by the Cauchy--Schwarz
inequality,
\begin{eqnarray*}
\bigl|l_{f_2}(f_1)\bigr|\leq\sqrt{\Var Z_1(f_1)
\Var Z_2(f_2)} =\Vert f_1\Vert_{\mathcal{H}(S_1)}
\Vert f_2\Vert_{\mathcal{H}(S_2)}.
\end{eqnarray*}
Thus, $l_{f_2}$ is a bounded linear functional on $\mathcal{H}(S_1)$
and by the Riesz representation theorem there is a bounded operator
$C_{12}\dvtx \mathcal{H}(S_2)\rightarrow\mathcal{H}(S_1)$ satisfying
%
\begin{equation}
\label{eqcompform} \Cov\bigl(Z_1(f_1),Z_2(f_2)
\bigr)= \langle f_1,C_{12}f_2\rangle
_{\mathcal{H}(S_1)}.
\end{equation}
There is also a bounded operator $C_{21}\dvtx  \mathcal{H}(S_1)\rightarrow
\mathcal{H}(S_2)$ with $C_{21}=C_{12}^*$, which satisfies
$\Cov(Z_1(f_1),Z_2(f_2))=\langle C_{21}f_1,f_2\rangle _{\mathcal{H}(S_2)}$.

%
\begin{proposition}
\label{proc12}
Under Assumption~\ref{aa1}, $\Vert C_{12}\Vert=\Vert C_{21}\Vert<1$.
\end{proposition}

\begin{pf} By the definition of the operator norm, we have
\begin{eqnarray*}
\Vert C_{12}\Vert^2=\sup_{f_2\in
\mathcal{H}(S_2),\Vert f_2\Vert_{\mathcal{H}(S_2)}=1}\Vert
C_{12}f_2\Vert^2_{\mathcal
{H}(S_1)}.
\end{eqnarray*}
An application of the Cauchy--Schwarz inequality produces
\begin{eqnarray*}
\bigl|\Cov\bigl(Z_1(f_1),Z_2(f_2)
\bigr)\bigr|&=& \bigl|\langle f_1,C_{12}f_2\rangle
_{\mathcal
{H}(S_1)}\bigr|
\\
&<& \sqrt{\Var Z_1(f_1)\Var Z_2(f_2)}
\\
&=& \Vert f_1\Vert_{\mathcal{H}(S_1)}\Vert f_2
\Vert_{\mathcal{H}(S_2)}
\end{eqnarray*}
with the strict inequality coming from Assumption~\ref{aa1}.
Now take $f_1=C_{12}f_2$.
\end{pf}

The operators $C_{12}$ and $S_{12}$ are, of course, related as we now
explain. For this purpose, define
\begin{eqnarray*}
\tilde{\mathcal{H}}(S_i)= \Biggl\{\tilde{f}_i\dvt
\tilde{f}_i=\sum_{j=1}^\infty
\tilde{f}_{ji}\phi_{ij}, \Vert\tilde{f}_i
\Vert^2_{\tilde{\mathcal{H}}(S_i)}=\sum_{j=1}^\infty
\lambda_{ij}\tilde{f}_{ij}^2 = \bigl\Vert
S_i^{1/2}\tilde{f}_i\bigr\Vert^2 <
\infty \Biggr\},\qquad  i = 1, 2.
\end{eqnarray*}
Then, $S_i$ is an isometric mapping from $\tilde{\mathcal{H}}(S_i)$
onto $\mathcal{H}(S_i)$; that is, $\tilde{\mathcal{H}}(S_i) =
S_i^{-1}\mathcal{H}(S_i)$. This leads us to
the following lemma.

%
\begin{lemma}
\label{lemkhatrilem}
$S_{12}$ is an operator from $\tilde{\mathcal{H}}(S_2)$
into $\mathcal{H}(S_1)$ with $\Vert S_{12}\Vert < 1$.
\end{lemma}

\begin{pf} For any $\tilde{f}_2 \in\tilde{\mathcal{H}}(S_2)$ and
$f_1 \in\mathcal{H}(S_1)$
\begin{eqnarray*}
\Cov\bigl(Z_1(f_1), Z_2(S_2
\tilde{f}_2)\bigr) & = & \sum_{i,j}f_{1i}
\tilde{f}_{2j}\langle\phi_{1i}, S_{12}
\phi_{2j} \rangle
\\
& = & \sum_{i,j}f_{1i}
\tilde{f}_{2j}\bigl\langle S_1^{1/2}
\phi_{1i}, S_1^{1/2}S_{12}
\phi_{2j} \bigr\rangle_{\mathcal{H}(S_1)}
\\
& = & \sum_{i,j}\lambda_{1i}f_{1i}
\tilde{f}_{2j}\langle\phi_{1i}, S_{12}
\phi_{2j} \rangle_{\mathcal{H}(S_1)}
\\
& = & \langle f_{1}, S_{12} \tilde{f}_2
\rangle_{\mathcal{H}(S_1)}.
\end{eqnarray*}
Now use the Cauchy--Schwarz inequality and
$\Vert S_2\tilde{f}_2\Vert_{\mathcal{H}(S_2)} =
\Vert\tilde{f}_2\Vert_{\tilde{\mathcal{H}}(S_2)}$.
\end{pf}

Lemma~\ref{lemkhatrilem} provides the means to characterize
$C_{12}$. Specifically, observe that
\begin{eqnarray*}
\Cov\bigl(Z_1(f_1), Z_2(S_2
\tilde{f}_2)\bigr) & = & \langle f_{1}, S_{12}
\tilde{f}_2 \rangle_{\mathcal{H}(S_1)}
\\
& = & \bigl\langle f_{1}, S_{12}S_2^{-1}S_2
\tilde{f}_2 \bigr\rangle_{\mathcal{H}(S_1)}
\\
& = & \langle f_{1}, C_{12}S_2
\tilde{f}_2 \rangle_{\mathcal{H}(S_1)}.
\end{eqnarray*}
In addition, the fact that $S_{12}$ is compact on $\mathcal{H}$ along
with an argument similar to that of Lemma~\ref{lemkhatrilem} reveals
that $C_{12}$ is the limit of a sequence of finite dimensional
operators. We summarize these findings as follows.

\begin{theorem}
\label{thsvd}
$C_{12} = S_{12}S_2^{-1} $ is a compact operator from
$\mathcal{H}(S_2)$ into $\mathcal{H}(S_1)$.
\end{theorem}

For $h\in\mathcal{H}_0$, define
$Qh=(f_1+C_{12}f_2,f_2+C_{21}f_1)$. It will be convenient to write
this in matrix form as
%
\begin{eqnarray} \label{eQhmatrix}
Qh=\lleft[ %
\begin{array} {c@{\quad}c}I& C_{12}
\\
C_{21}& I \end{array} %
\rright]\lleft[ %
\begin{array} {c} f_{1}
\\
f_2 \end{array} %
\rright]
\end{eqnarray}
with the convention that the resulting vector is viewed as an element
of $\mathcal{H}_0$. Observe that
\begin{eqnarray*}
\Cov\bigl(Z(h),Z\bigl(h'\bigr)\bigr) = \bigl\langle
h,Qh'\bigr\rangle _0.
\end{eqnarray*}
This leads to the following proposition.

%
\begin{proposition}
\label{proQinv}
$Q\dvtx  \mathcal{H}_0\rightarrow\mathcal{H}_0$ is invertible with inverse
defined by
%
\begin{eqnarray}\label{eQinv}
Q^{-1}(h)=\bigl(C_{11.2}^{-1}f_1-C_{12}C_{22.1}^{-1}f_2,C_{22.1}^{-1}f_2-C_{21}C_{11.2}^{-1}f_1
\bigr),
\end{eqnarray}
where $h=(f_1,f_2)\in\mathcal{H}_0$ and
$C_{ii.k}=I-C_{ik}C_{ki}=(I-C_{ik}C_{ki})^*$
for $i,k=1,2,i \neq k$.
\end{proposition}

Analogous to \eqref{eQhmatrix}, \eqref{eQinv} will also be expressed as
\begin{eqnarray*}
Q^{-1}h=\lleft[ %
\begin{array} {c@{\quad}c}C_{11.2}^{-1}&
-C_{12}C_{22.1}^{-1}
\\\noalign{\vspace*{2pt}}
-C_{21}C_{11.2}^{-1}& C_{22.1}^{-1}
\end{array} %
\rright]\lleft[ %
\begin{array} {c}
f_{1}
\\
f_2 \end{array} %
\rright].
\end{eqnarray*}

\begin{pf}
The form of the inverse as stated in \eqref{eQinv} follows directly
once we have shown all the relevant inverse operators exist. Thus, let
us concentrate on the latter task.

We can write $Q=I-T$ with
\begin{eqnarray*}
Th=(-C_{12}f_2,-C_{21}f_1)=-\lleft[
\begin{array} {c@{\quad}c}0& C_{12}
\\
C_{21}& 0 \end{array} %
\rright] \lleft[ %
\begin{array} {c} f_{1}
\\
f_2 \end{array} %
\rright].
\end{eqnarray*}
Then
\begin{eqnarray*}
\Vert Th\Vert^2_0&=&\Vert C_{12}f_2
\Vert^2_{\mathcal
{H}(S_1)}+\Vert C_{21}f_1
\Vert^2_{\mathcal{H}(S_2)}
\\
&\leq& \Vert C_{12}\Vert^2\Vert|f_2
\Vert^2_{\mathcal
{H}(S_2)}+\Vert C_{21}\Vert^2\Vert
f_1\Vert^2_{\mathcal{H}(S_1)}
\\
&=& \Vert C_{12}\Vert^2\bigl[\Vert f_1
\Vert^2_{\mathcal{H}(S_1)}+\Vert f_2\Vert^2_{\mathcal
{H}(S_2)}
\bigr]
\\
&=& \Vert C_{12}\Vert^2\Vert h\Vert_0^2
\\
&<& \Vert h\Vert_0^2
\end{eqnarray*}
by Proposition~\ref{proc12}. Theorem~4.40 of Rynne and Youngson \cite
{randy} now has the
consequence that $I-T=Q$ is invertible.

To complete the proof, we need to show that $C_{11.2}$ and $C_{22.1}$
are invertible. This again follows from Theorem~4.40 of Rynne and Youngson \cite{randy}
because $C_{11.2}=I-C_{12}C_{21}$ with $\Vert C_{21}\Vert=\Vert
C_{12}\Vert<1$ from
Proposition~\ref{proc12}.
\end{pf}

Now define
\begin{eqnarray*}
\mathcal{H}(Q)=\lleft\{h\dvt  h=Q\lleft[ %
\begin{array} {c}
f_{1}
\\
f_2 \end{array} %
\rright], f_i\in
\mathcal{H}(S_i),i=1,2,\Vert h\Vert^2_{\mathcal{H}(Q)}=
\bigl\Vert Q^{-1/2}h\bigr\Vert_0^2 < \infty \rright\}.
\end{eqnarray*}
The next proposition follows immediately from this definition.

%
\begin{proposition}{
\label{proHQ}}
$\mathcal{H}(Q)$ is congruent to
\begin{eqnarray*}
L_Z^2= \bigl\{Z(h)\dvt  h\in \mathcal{H}_0,
\bigl\Vert Z(h)\bigr\Vert^2_{L^2_Z} := \Var\bigl(Z(h)\bigr) < \infty
\bigr\}
\end{eqnarray*}
under the mapping $\Psi(h) = Z(Q^{-1}h)$.
\end{proposition}

With Proposition~\ref{proHQ} in hand we can now give our formulation
of CCA. Specially, we seek elements $f_i \in\mathcal{H}(S_i)$ of
unit norm that maximize $|\Cov(Z_1(f_1),Z_2(f_2))|$. But
\begin{eqnarray*}
\Cov\bigl(Z_1(f_1),Z_2(f_2)
\bigr)=\Cov\bigl(Z(f_1,0),Z(0,f_2)\bigr) = \left\langle Q
\lleft[ %
\begin{array} {c} f_{1}
\\
0 \end{array} %
\rright],Q\lleft[ %
\begin{array} {c} 0
\\
f_2 \end{array} %
\rright] \right\rangle_{\mathcal{H}(Q)}
\end{eqnarray*}
which leads to the conclusion that it is equivalent to find $f_i\in
\mathcal{H}(S_i)$ to maximize the right-hand side of this last
expression.

The analysis from this point is driven by the results of Sunder \cite{sunder}
as described in Section~\ref{sec4}. For that purpose, we decompose
$\mathcal{H}(Q)$ into a sum of the closed subspaces
$M_1$ and $M_2$ with
\begin{eqnarray*}
M_1&=&\lleft\{h\in\mathcal{H}(Q)\dvt  h=Q\lleft[ %
\begin{array} {c} f_{1}
\\
0 \end{array} %
\rright]:=(f_1,C_{21}f_1),f_1
\in\mathcal{H}(S_1) \rright\},
\\
M_2&=&\lleft\{h \in\mathcal{H}(Q)\dvt  h=Q\lleft[ %
\begin{array} {c} 0
\\
f_2 \end{array} %
\rright]:=(C_{12}f_2,f_2),f_2
\in \mathcal{H}(S_2) \rright\}.
\end{eqnarray*}
Regarding $M_1$ and $M_2$, we have the following result.

%
\begin{proposition}
\label{proHQsum}
$\mathcal{H}(Q)=M_1+M_2$ with ``$+$'' indicating an algebraic direct sum.
\end{proposition}

\begin{pf} Clearly any element of $\mathcal{H}_0$ can be written
as the sum of elements in $M_1$ and $M_2$. We therefore need only
show that $M_1 \cap M_2= \{0 \}$. Thus, suppose there exist
$f_i\in\mathcal{H}(S_i), i = 1, 2$, such that $(f_1,C_{21}f_1) =
(C_{12}f_2,f_2)$.
Then
\begin{eqnarray*}
\Var\bigl(Z_1(f_1)\bigr)=\langle f_1,f_1
\rangle _{\mathcal
{H}(S_1)}=\langle f_1,C_{12}f_2
\rangle _{\mathcal{H}(S_1)}
\end{eqnarray*}
and
\begin{eqnarray*}
\Var\bigl(Z_2(f_2)\bigr)=\langle f_2,f_2
\rangle _{\mathcal
{H}(S_2)}=\langle f_2,C_{21}f_1
\rangle _{\mathcal{H}(S_2)} =\langle C_{12}f_2,f_1
\rangle _{\mathcal{H}(S_1)}.
\end{eqnarray*}
But, these relations have the consequence that
$|\Corr(Z_1(f_1),Z_2(f_2))|=1$ which contradicts Assumption~\ref{aa1}.
\end{pf}

To relate Proposition~\ref{proHQsum} to Sunder's scheme in the
\hyperref[sec5]{Appendix}, let $L_1=M_1$ and $L_2=M_2\cap M_1^\bot$ in Theorem~\ref{thsunder}. Then, for $h_1=Q
\bigl[{f_{1} \atop 0}\bigr]\in M_1$ and $h_2=Q\bigl[{0\atop f_2}\bigr]\in M_2$, the first canonical correlation
satisfies
\begin{eqnarray*}
\rho&=&\mathop{\sup_{h_1\in M_1,h_2\in M_2}}_{\Vert h_i\Vert_{\mathcal{H}(Q)}=
1, i= 1,2}\bigl|\langle h_1,h_2
\rangle_{\mathcal{H}(Q)}\bigr| = \mathop{\sup_{h_1\in
L_1,\tilde{h}_2\in L_2}}_{\Vert h_1\Vert_{\mathcal{H}(Q)}= 1,
\Vert\tilde{h}_2 + B\tilde{h}_2\Vert_{\mathcal{H}(Q)}= 1}\bigl|\langle h_1,B
\tilde{h}_2 \rangle_{\mathcal{H}(Q)}\bigr|
\\
&\le&\mathop{\sup_{\tilde{h}_2\in L_2}}_{\Vert\tilde{h}_2 +
B\tilde{h}_2\Vert_{\mathcal{H}(Q)}= 1}\Vert B\tilde{h}_2\Vert
_{\mathcal{H}(Q)}
\end{eqnarray*}
for $B=P_{L_1|M_2}(P_{L_2|M_2})^{-1}$. Taking
$h_1=B\tilde{h}_2/\Vert B\tilde{h}_2\Vert_{\mathcal{H}(Q)}$, we see
that the
bound is attainable and holds with equality. Thus, we have shown that
$\rho$ is obtained by maximizing $\Vert B\tilde{h}_2\Vert_{\mathcal{H}(Q)}$
subject to
\begin{eqnarray*}
\Vert B\tilde{h}_2+\tilde{h}_2\Vert_{\mathcal{H}(Q)}^2
=\bigl\langle \tilde{h}_2,\bigl(I+B^*B\bigr)\tilde{h}_2
\bigr\rangle _{\mathcal{H}(Q)} = 1.
\end{eqnarray*}

The operator $B^{*}B$ is compact as a result of Theorem~\ref{thsvd}
and Theorem~\ref{thccachar} below. In addition, $I+B^*B$ is
self-adjoint, positive, invertible and has a self-adjoint square-root
$(I+B^*B)^{1/2}$. We can therefore work with
$\tilde{h}_2'=(I+B^*B)^{1/2}\tilde{h}_2$ and maximize
\begin{eqnarray*}
\Vert B\tilde{h}_2\Vert_{\mathcal{H}(Q)}=\bigl\Vert B\bigl(I+B^*B
\bigr)^{-1/2}\tilde {h}_2'\bigr\Vert_{\mathcal{H}(Q)}
\end{eqnarray*}
subject to $\tilde{h}_2' \in L_2$ and
$\Vert\tilde{h}_2'\Vert_{\mathcal{H}(Q)}^2=1$. The maximizer is the
eigenvector for the largest eigenvalue of
$(I+B^*B)^{-1/2}B^*B(I+B^*B)^{-1/2}$. Some algebra reveals that the
resulting eigenvalue problem is equivalent to finding a vector
$\tilde{h}_2 \in L_2$ with $\Vert\tilde{h}_2\Vert_{\mathcal
{H}(Q)}^2=1$ such
that
%
\begin{eqnarray}
\label{e2} B^*B\tilde{h}_2=\alpha^2
\tilde{h}_2
\end{eqnarray}
in which case $\rho=\alpha/\sqrt{1+\alpha^2}$.

Now suppose that $\tilde{h}_2 \in L_2$ is any vector that satisfies
\eqref{e2}. Its $M_1$ component is $B\tilde{h}_2$ and its $M_2$
component is $B\tilde{h}_2+\tilde{h}_2$. These correspond to the
canonical variables $\Psi (B\tilde{h}_2/\alpha )$ and
$\Psi
((\tilde{h}_2+B\tilde{h}_2)/\sqrt{1+\alpha^2} )$ of the
$Z_1$ and $Z_2$ spaces, respectively.

In combination Corollaries \ref{corBh} and \ref{corB*} from the
\hyperref[sec5]{Appendix} give us the desired characterization for $B^*B$: namely,

\begin{theorem}
\label{thccachar}
For $h=(0,\tilde{f}_2) \in L_2, B^*B(0,\tilde{f_2}) =
(0,C_{21}C_{12}C_{22.1}^{-1}\tilde{f_2})$.
\end{theorem}

An application of Proposition~\ref{proPh} from the \hyperref[sec5]{Appendix} now
reveals that the conclusion of Theorem~\ref{thccachar} can be
restated as $B^*B(0,\tilde{f_2}) =
(0,C_{21}C_{12}f_2)$ for some $f_2 \in\mathcal{H}(S_2)$ and the
eigenvalue problem (\ref{e2}) is equivalent to $C_{21}C_{12}f_2 =
\alpha^2 C_{22.1}f_2$ or
\begin{eqnarray*}
C_{21}C_{12}f_2 = \rho^2f_2.
\end{eqnarray*}
By interchanging the roles of $M_1$ and $M_2$ it follows that the
optimal choice for $f_1$ is the eigenvector corresponding to the same
eigenvalue $\rho^2$ of $C_{12}C_{21}$. Thus, $\rho$ is the largest
singular value of $C_{21}$, $f_1, f_2$ are its right and left hand
singular functions and $Z_1(f_1), Z_2(f_2)$ are the corresponding
canonical variables. More generally, a similar analysis reveals that
the collection of all such singular values gives rise to a sequence of
canonical correlations that correspond to canonical variable pairs
with maximum possible correlation subject to being uncorrelated with
previous pairs in the sequence.

We conclude this section with examples that illustrate some of the
features of our CCA formulation.

%
\begin{example}
Suppose that $S_1$ and $S_2$ are full-rank, finite-dimensional
matrices. Then, $C_{12}=S_{12}S_2^{-1}$ and $C_{21} = S_{21}S_1^{-1}$
so that finding eigenvalues and eigenvectors for $C_{21}C_{12}$ is
equivalent to the singular value decomposition of
$S_1^{-1/2}S_{12}S_2^{-1/2}$ which, in turn, is equivalent to
Hotelling's classic solution for the finite dimensional case as
established in Kshirsagar \cite{ksh}.
\end{example}

%
\begin{example}
Functional data analysis generally focuses on the case where the $X_i$
are random element of $L^2[0, 1]$; that is, the set of square integrable
function on the interval $[0, 1]$. One assumes the $X_i$ admit point-wise
representations as the continuous time stochastic processes $\{X_i(t,
\omega) \dvt  t \in[0, 1], \omega\in\Omega\}, i = 1, 2$. Inference is then
based on the linear combinations described in (\ref{eql2subx}).

The (assumed continuous) process covariance kernels are
\begin{eqnarray*}
K_i\bigl(t,t'\bigr)=\Cov\bigl(X_i(t),X_i
\bigl(t'\bigr)\bigr)= \sum_{j=1}^\infty
\lambda_{ij}\phi_{ij}(t)\phi_{ij}
\bigl(t'\bigr)
\end{eqnarray*}
with the $(\lambda_{ij}, \phi_{ij}), j = 1, \ldots, i = 1, 2$, being
the eigenvalues and eigenvectors of the $L^2[0, 1]$ integral operators
defined by
\begin{eqnarray*}
(S_{i}f) (t) = \int_0^1f(s)K_i(t,
s)\,\mathrm{d}s.
\end{eqnarray*}
The RKHS that is congruent to $L^2_{X_i}$ is $\mathcal{H}(S_i)$.

In the case of two processes, we also have the cross-covariance kernels
\begin{eqnarray*}
K_{12}(t_1,t_2) & = &\Cov\bigl(X_1(t_1),X_2(t_2)
\bigr)
\\
&=&\Cov\bigl(X_2(t_2),X_1(t_1)
\bigr)
\\
&=& K_{21}(t_2,t_1) .
\end{eqnarray*}
From Eubank and Hsing \cite{eandh}, we know that $K_{12}(\cdot
,t_2) \in
\mathcal{H}(S_1)$, and $K_{12}(t_1,\cdot)\in\mathcal{H}(S_2)$; so,
if $f_i =\sum_{j = 1}^{\infty}\lambda_{ij}f_{ij}\phi_{ij} \in
\mathcal{H}(S_i)$,
\begin{eqnarray*}
(R_{12}f_2) (t)=\bigl\langle K_{12}(t,
\cdot),f_2(\cdot)\bigr\rangle _{\mathcal{H}(S_2)}
\end{eqnarray*}
defines a bounded operator from $\mathcal{H}(S_2)$ into
$\mathcal{H}(S_1)$ with the property that
\begin{eqnarray*}
\Cov \bigl(Z_1(f_1), Z_2(f_2)
\bigr) &=&\sum_k\sum_jf_{1j}f_{2k}
\int_0^1K_{12}(s, t)
\phi_{1j}(s)\phi _{2k}(t)\,\mathrm{d}s\,\mathrm{d}t
\\
&=& \langle f_1,R_{12}f_2\rangle_{\mathcal{H}(S_1)}.
\end{eqnarray*}
Therefore, $R_{12}=C_{12}$ and our CCA formulation coincides with that
in Eubank and Hsing \cite{eandh}.
\end{example}

%
\begin{example}
\label{examcomp}
The developments in this section suggest a new approach to estimation
in the functional CCA setting of the previous example. The idea stems from
(\ref{eqcompform}) which has the consequence that
%
\begin{equation}
\label{eqcompform2} \Cov\bigl(Z_1(\phi_{1i}),Z_2(
\phi_{2j})\bigr)= \langle \phi_{1i},C_{12}
\phi_{2j}\rangle _{\mathcal{H}(S_1)}.
\end{equation}
It follows from Hansen \cite{hansen} that a singular
value decomposition of
%
\begin{equation}
\label{eqasubm} A_m = \bigl\{ \langle \phi_{1i},C_{12}
\phi_{2j}\rangle _{\mathcal{H}(S_1)} \bigr\}_{i, j = 1:m}
\end{equation}
for some finite integer $m$ will produce singular values that
approximate the singular values for the operator $C_{12}$ and that the
singular vectors provide coefficients for linear combinations of the
$\phi_{ij}$ that approximate its singular functions. The only question
is how to estimate the inner products in (\ref{eqasubm}). The answer
is revealed by examining the left hand of (\ref{eqcompform2}). The
realized values of the $Z_i(\phi_{ij}), j = 1, \ldots, m$ can be
estimated directly using the scores one obtains from a principal
components analysis of functional data. Thus, their sample covariance
matrix provides an obvious choice for an estimator of
(\ref{eqasubm}).

Suppose we have observed sample path pairs $ ( x_{1j}(\cdot),
x_{2j}(\cdot)  ), j = 1, \ldots, n$. The resulting estimation
algorithm can then be summarized as follows.
\begin{enumerate}[2.]
\item[1.] Carry out a principal components analysis of the $x_{ij}, j = 1,
\ldots, n$ to obtain the estimated eigenfunctions $\hat{\phi}_{ij},
j =
1,\ldots, m$ and $n \times m$ score matrices
\begin{eqnarray*}
W_i = \bigl\{ \bigl\langle\hat{\phi}_{ij},
x_{ik}(\cdot) \bigr\rangle \bigr\}_{k = 1:n, j = 1:m}
\end{eqnarray*}
for $i = 1, 2$. Let $\hat{A}_m$ be the $m \times m$ sample
cross covariance matrix obtained from $W_1$ and $W_2$.
\item[2.] If $\hat{A}_m = UDV^T$ for $U =  [u_1, \ldots, u_m
], V
=  [ v_1, \ldots, v_m  ]$ and $D = \diag(d_1, \ldots,
d_m)$ is the singular value decomposition of $\hat{A}_m$, the
$i$th canonical correlation is estimated by $d_i$ and the
corresponding canonical weight functions by $u_i^T[\hat{\phi}_{21},
\ldots, \hat{\phi}_{2m}]$ and $v_i^T[\hat{\phi}_{11}, \ldots,
\hat{\phi}_{1m}]$.
\end{enumerate}

A simple numerical example will be used to illustrate this estimation
scheme. The setting is that of Eubank and Hsing \cite{eandh} where
the two processes are
\begin{eqnarray*}
X_1(t) &=& \sum_{j=1}^{20}
j^{-1/2} Z_{1j} \sqrt{2}\sin(j\pi t),
\\
X_2(t) &=& (Z_{11} + Z_{21}) \sin(\pi t)+\sum
_{j=2}^{20} j^{-1/2} Z_{2j}
\sqrt{2}\sin(j\pi s),
\end{eqnarray*}
for $t\in[0,1]$ and the $Z_{ij}$ i.i.d. standard normal
random variables. In this instance, there is only one nonzero canonical
correlation: namely, $\rho_1 = 1/\sqrt{2} \doteq0.707$.

We sampled $n$ process pairs at 100 equally spaced points and
conducted principal components analysis on the resulting data using
the function \texttt{\textmd{pda.fd}} from the fda package in R
retaining 9 components (or harmonics) for both processes. This basic
experiment was then replicated 100 times. For samples of size $n =
250$, the observed means (standard deviations) of the first two sample
canonical correlations were 0.7248 (0.0818) and 0.0777 (0.0122),
respectively. For samples of size $n = 500$, the means (standard
deviations) were 0.7147 (0.0591) and 0.055 (0.0095).

This rather crude implementation suffices for the present expository
purposes. However, for use in practice one should at least employ
consistent estimators for the eigenfunctions such as those studied in
Yao \textit{et al.} \cite{yaoetal} and Hall \textit{et al.} \cite
{halletal}.
\end{example}

\section{PCCA}
\label{sec3}
A similar approach to that of the previous section can be used to
address the PCCA setting. There are now three $\mathcal{H}$-valued
random variables $X_i, i = 1, 2, 3$, with associated covariance
operators $S_i$, $i=1,2,3$. As in Section~\ref{sec2}, we can also
define the cross-covariance operators $S_{12}$, $S_{13}$, $S_{23}$ and
their adjoints.

For $i=1,2,3$, the Hilbert spaces $L_{Z_i}^2$ spanned by the process
$Z_i(f_i)$ indexed by their congruent Hilbert spaces
$\mathcal{H}(S_i)$ are defined as in \eqref{eLZ} and \eqref{eHS}.
Hence, by the Riesz representation theorem, there are bounded
operators $C_{ij}\dvtx  \mathcal{H}(S_j) \rightarrow\mathcal{H}(S_i)$
satisfying
\begin{eqnarray*}
\Cov\bigl(Z_i(f_i),Z_j(f_j)
\bigr)=\langle f_i, C_{ij}f_j
\rangle_{\mathcal{H}(S_i)}
\end{eqnarray*}
for $i,j=1,2,3$ and $i \neq j$. Also, we have that $C_{ij}=C_{ji}^*$.

We now construct the new Hilbert space
\begin{eqnarray*}
\mathcal{H}_0=\Biggl\{h=(f_1,f_2,f_3)\dvt
f_i\in\mathcal{H}(S_i), i=1,2,3, \|h
\|_0^2= \sum_{i = 1}^3
\|f_i\|_{\mathcal{H}(S_i)}^2 < \infty\Biggr\}.
\end{eqnarray*}
Then, our corresponding $\mathcal{H}_0$ indexed process is
$Z(h)=\sum_{i = 1}^3Z_i(f_i)$.

As in the previous section we need to rule out the case where
perfect prediction is possible. For this purpose, we require that
Assumption~\ref{aa1} holds for both of the process pairs
$Z_1$, $Z_2$ and $Z_1$, $Z_3$ as well as
the following.

%
\begin{assumption}
\label{aa2}
There exist no $f_2 \in\mathcal{H}(S_2)$ or $f_3 \in
\mathcal{H}(S_3)$ such that
\begin{eqnarray*}
\bigl|\Corr\bigl(Z_2(f_2)-P_{Z_1}Z_2(f_2),Z_3(f_3)-P_{Z_1}Z_3(f_3)
\bigr)\bigr|=1.
\end{eqnarray*}
\end{assumption}

For $h\in\mathcal{H}_0$, define
\begin{eqnarray*}
Qh=(f_1+C_{12}f_2+C_{13}f_3,
C_{21}f_1+f_2+C_{23}f_3,C_{31}f_1+C_{32}f_2+f_3)
\end{eqnarray*}
which we will express in the matrix form
\begin{eqnarray*}
Qh=\lleft[ %
\begin{array} {c@{\quad}c@{\quad}c}I &C_{12} & C_{13}
\\
C_{21} & I & C_{23}
\\
C_{31} & C_{32} & I \end{array} %
\rright]\lleft[
\begin{array} {c}f_1
\\
f_2
\\
f_3 \end{array} %
\rright].
\end{eqnarray*}
We then see that
\begin{eqnarray*}
\Cov\bigl(Z(h),Z\bigl(h'\bigr)\bigr)=\bigl\langle
h,Qh'\bigr\rangle_0.
\end{eqnarray*}

Our next result gives the three process parallel of Proposition~\ref{proQinv}.

%
\begin{proposition}
Let $E =[ {C_{12} \enskip
C_{13}
}]$,
$F= \bigl[{C_{21}\atop C_{31}}
\bigr]$,
$D=\bigl[{
I \atop C_{23}}\enskip{ C_{32} \atop
I}\bigr]$ and
$G= D^{1/2}(I- V )D^{1/2}$
with
%
\begin{equation}
\label{eQC} V = \lleft[ %
\begin{array} {c@{\quad}c}0 &
-C_{22.1}^{-1/2}(C_{23}-C_{21}C_{13})C_{33.1}^{-1/2}
\\\noalign{\vspace*{2pt}}
-C_{33.1}^{-1/2}(C_{32}-C_{31}C_{12})C_{22.1}^{-1/2}
& 0 \end{array} %
\rright].
\end{equation}
Then,
%
\begin{eqnarray}
Q^{-1}=\lleft[ %
\begin{array} {c@{\quad}c}I+EG^{-1}F
& -EG^{-1}
\\
G^{-1}F & G^{-1} \end{array} %
\rright].
\label{eQinv2}
\end{eqnarray}
\end{proposition}

From Proposition~\ref{proQinv}, we know that $C_{22.1}$ and
$C_{33.1}$ are invertible. The result will therefore follow if we can
show that the norm of $V$ in \eqref{eQC} is strictly less than
unity. This is a consequence of the next two lemmas and Theorem~4.40
of Rynne and Youngson \cite{randy}.

%
\begin{lemma}
\label{lemprojlem}
The projection of $Z_2(f_2)$ onto $L^2_{Z_1}$ is $Z_1(C_{12}f_2)$
and the projection of $Z_3(f_3)$ onto $L^2_{Z_1}$ is $Z_1(C_{13}f_3)$.
\end{lemma}

\begin{pf} If $P_{Z_1}Z_2(f_2)$ denotes the projection, it must satisfy
\begin{eqnarray*}
\Cov\bigl(Z_1(f_1),P_{Z_1}Z_2(f_2)
\bigr)=\Cov\bigl(Z_1(f_1),Z_2(f_2)
\bigr)
\end{eqnarray*}
for every $f_1 \in\mathcal{H}(S_1)$. Since there is some $\check
{f}_1\in
\mathcal{H}(S_1)$ such that $P_{Z_1}Z_2(f_2)=Z_1(\check{f}_1)$,
\begin{eqnarray*}
\Cov\bigl(Z_1(f_1),Z_2(f_2)
\bigr)&=&\langle f_1,C_{12}f_2 \rangle
_{\mathcal{H}(S_1)}
\\
&=& \Cov\bigl(Z_1(f_1),Z_1(
\check{f}_1)\bigr)
\\
&=& \langle f_1, \check{f}_1 \rangle_{\mathcal{H}(S_1)}.
\end{eqnarray*}
Therefore, $\check{f}_1 = C_{12}f_2$. The second half of the lemma is proved
similarly.
\end{pf}

\begin{lemma}
\label{lemnorm}
$\|C_{22.1}^{-1/2}(C_{23}-C_{21}C_{13})C_{33.1}^{-1/2}\|_{\mathcal{H}(S_2)}<1$.
\end{lemma}

\begin{pf} First, observe that by Lemma~\ref{lemprojlem}
and Assumption~\ref{aa2}
\begin{eqnarray*}
&&\bigl|\Cov\bigl(Z_2(f_2)-Z_1(C_{12}f_2),Z_3(f_3)-Z_1(C_{13}f_3)
\bigr)\bigr|
\\
&&\quad = \bigl|\langle f_2,C_{23}f_3
\rangle_{\mathcal{H}(S_2)}-\langle f_2, C_{21}C_{13}f_3
\rangle_{\mathcal{H}(S_2)}\bigr|
\\
&&\quad < \bigl(\Var\bigl(Z_2(f_2)-Z_1(C_{12}f_2)
\bigr)\bigr)^{1/2}\bigl(\Var \bigl(Z_3(f_3)-Z_1(C_{13}f_3)
\bigr)\bigr)^{1/2}
\\
&&\quad  = \langle f_2, C_{22.1}f_2
\rangle_{\mathcal{H}(S_2)}^{1/2}\langle f_3, C_{33.1}f_3
\rangle_{\mathcal{H}(S_3)}^{1/2}
\\
&&\quad  = \bigl\|C_{22.1}^{1/2}f_2\bigr\|_{\mathcal{H}(S_2)}
\bigl\|C_{33.1}^{1/2}f_3\bigr\| _{\mathcal
{H}(S_3)}.
\end{eqnarray*}
Now, let $\tilde{f}_2=C_{22.1}^{1/2}f_2$ and
$\tilde{f}_3=C_{33.1}^{1/2}f_3$ to obtain
\begin{eqnarray*}
\bigl\langle\tilde{f}_2, C_{22.1}^{-1/2}(C_{23}-C_{21}C{13})C_{33.1}^{-1/2}
\tilde{f}_3 \bigr\rangle_{\mathcal{H}(S_2)} < \|\tilde{f}_2
\|_{\mathcal{H}(S_2)}\|\tilde{f}_3\|_{\mathcal{H}(S_3)}.
\end{eqnarray*}
Finally, taking
$\tilde{f}_2=C_{22.1}^{-1/2}(C_{23}-C_{21}C_{13})C_{33.1}^{-1/2}\tilde{f}_3$
completes the proof.
\end{pf}

Now define
\begin{eqnarray*}
\mathcal{H}(Q)=\lleft\{h\dvt  h=Q\lleft[ %
\begin{array}
{c}f_1
\\
f_2
\\
f_3 \end{array} %
\rright], f_i \in
\mathcal{H}(S_i), i=1,2,3, \|h\|^2_{\mathcal{H}(Q)}=
\bigl\|Q^{-1/2}h\bigr\|_0^2 < \infty \rright\}.
\end{eqnarray*}
Then, as in Proposition~\ref{proHQ}, we have

%
\begin{proposition}
$\mathcal{H}(Q)$ is congruent to
\begin{eqnarray*}
L_Z^2=\bigl\{Z(h)\dvt  h\in \mathcal{H}_0,
\bigl\|Z(h)\bigr\|^2_{L_Z^2}=\Var\bigl(Z(h)\bigr) < \infty\bigr\}
\end{eqnarray*}
under the mapping $\Psi(h) = Z(Q^{-1}h)$.
\end{proposition}

For the PCCA formulation, we wish to find $f_2\in\mathcal{H}(S_2)$
and $f_3 \in\mathcal{H}(S_3)$ to maximize
\begin{eqnarray*}
\bigl| \Cov\bigl(Z_2(f_2)-Z_1(C_{12}f_2),Z_3(f_3)-Z_1(C_{13}f_3)
\bigr)\bigr|.
\end{eqnarray*}
Since
\begin{eqnarray*}
&&\Cov\bigl(Z_2(f_2)-Z_1(C_{12}f_2),Z_3(f_3)-Z_1(C_{13}f_3)
\bigr)
\\
&&\quad =\Cov\bigl(Z(-C_{12}f_2,f_2,0),Z(-C_{13}f_3,0,f_3)
\bigr),
\end{eqnarray*}
it suffices to find $f_2\in\mathcal{H}(S_2)$ and $f_3 \in
\mathcal{H}(S_3)$ to maximize
\begin{eqnarray*}
\left\vert \left\langle Q\lleft[ %
\begin{array}
{c}-C_{12}f_2
\\
f_2
\\
0 \end{array} %
\rright],Q\lleft[ %
\begin{array}
{c}-C_{13}f_3
\\
0
\\
f_3 \end{array} %
\rright] \right\rangle_{\mathcal{H}(Q)}\right\vert .
\end{eqnarray*}

Again, we apply the results of Sunder described in Section~\ref{sec4}.
For this purpose, write $\mathcal{H}(Q)=M_1+M_2+M_3$
with
\begin{eqnarray*}
M_1=\lleft\{h\in\mathcal{H}(Q)\dvt  h=Q\lleft[ %
\begin{array} {c}f_1
\\
0
\\
0 \end{array} %
\rright]:=(f_1,C_{21}f_1,C_{31}f_1)
\rright\},
\end{eqnarray*}
\begin{eqnarray*}
M_2=\lleft\{h\in\mathcal{H}(Q)\dvt  h=Q\lleft[ %
\begin{array} {c}0
\\
f_2
\\
0 \end{array} %
\rright]:=(C_{12}f_2,f_2,C_{32}f_2)
\rright\}
\end{eqnarray*}
and
\begin{eqnarray*}
M_3=\lleft\{h\in\mathcal{H}(Q)\dvt  h=Q\lleft[ %
\begin{array} {c}0
\\
0
\\
f_3 \end{array} %
\rright]:=(C_{13}f_3,C_{23}f_3,f_3)
\rright\}.
\end{eqnarray*}
An argument similar to that for Proposition~\ref{proHQsum} produces
the following proposition.

%
\begin{proposition}
$\mathcal{H}(Q)=M_1+M_2+M_3$ with ``$+$'' indicating an algebraic
direct sum.
\end{proposition}

Now let $L_1=M_1$, $L_2=M_2\cap M_1^\bot, L_3=M_3 \cap M_2^\bot\cap
M_1^\bot$ and take
\begin{eqnarray*}
\hat{h}_2=Q\lleft[ %
\begin{array} {c}-C_{12}f_2
\\
f_2
\\
0 \end{array} %
\rright] \in M_2-P_{L_1}M_2
\quad \mbox{and}\quad  \hat{h}_3=Q\lleft[ %
\begin{array}
{c}-C_{13}f_3
\\
0
\\
f_3 \end{array} %
\rright] \in M_3-P_{L_1}M_3
\end{eqnarray*}
with $\|\hat{h}_i\|_{\mathcal{H}(Q)}=1$, $i=2,3$. Then, arguing as in
the previous section we see that the first partial canonical
correlation can be characterized as
\begin{eqnarray*}
\rho&=& \mathop{\sup_{\hat{h}_2\in M_2-P_{L_1}M_2, \hat{h}_3\in
M_3-P_{L_1}M_3}}_{\Vert\hat{h}_i\Vert_{\mathcal{H}(Q)} = 1, i
= 2,
3}\bigl\llvert \langle\hat{h}_2,
\hat{h}_3 \rangle_{\mathcal{H}(Q)} \bigr\rrvert
\\
&=& \sup_{\tilde{h}_3\in L_3,\Vert\tilde{h}_3 + B
\tilde{h}_3\Vert_{\mathcal{H}(Q)}| = 1 } \|B\tilde{h}_3\|_{\mathcal{H}(Q)}
\end{eqnarray*}
for $B=P_{L_2|M_3}(P_{L_3|M_3})^{-1}$. The bound is attained by taking
$\hat{h}_2=B\tilde{h}_3/\|B\tilde{h}_3\|_{\mathcal{H}(Q)}$ in which
case the first partial canonical correlation is $\alpha/\sqrt{1 +
\alpha^2}$ with $\alpha^2$ the largest eigenvalue of $B^*B$. If
$\tilde{h}_3$ is an eigenvector corresponding to $\alpha^2$, the
partial canonical variable for the $Z_2$ space is $\Psi
(B\tilde{h}_3/\alpha )$ and the partial canonical variable
for the $Z_3$ space is $\Psi
((\tilde{h}_3+B\tilde{h}_3)/\sqrt{1+\alpha^2} )$.

Now, through Corollaries \ref{proPBh} and \ref{proPB*h}, we finally obtain

%
\begin{theorem}
\label{thmpccabb}
For $h=(0,0,\tilde{f}_3) \in L_3$,
\begin{eqnarray*}
B^*Bh= \bigl(0,0,(C_{32}-C_{31}C_{12})C_{22.1}^{-1}(C_{23}-C_{21}C_{13})C_0^{-1}
\tilde {f}_3\bigr).
\end{eqnarray*}
\end{theorem}

This result in combination with Corollary~\ref{corl3form} reveals
that partial canonical correlations are the singular values of the
operator $C_{33.1}^{-1/2}(C_{32} - C_{31}C_{12})C_{22.1}^{-1/2}$.

%
\begin{example}
The basic computational algorithm from Example~\ref{examcomp} can be
adapted for computing sample partial canonical correlations. One now
carries out principal components analysis of the data from all three
processes and then regresses the scores for the $X_2, X_3$ process
data onto the scores from the $X_1$ sample paths. The Example~\ref{examcomp} computational scheme is then applied to the residuals
from the two regression analyses.

To illustrate the idea, consider again the two processes from Example~\ref{examcomp}. Sample paths were generated as before except
that in each instance we subtracted a term $\beta Z \cos( \pi s)$ with
$Z$ a standard normal random variable and $\beta$ equal to 1 for the
$X_1$ process and 2 for the $X_2$ process. The only nonzero partial canonical
correlation in this case is again $1/\sqrt{2}$. The first two partial
canonical correlations obtained from an empirical experiment using
the same parameters as in Example~\ref{examcomp} had means (standard
deviations) of 0.7107 (0.0875) and 0.0818 (0.0157) for samples of size
250 and 0.7141 (0.0599) and 0.0553 (0.0089) for samples of size 500.
\end{example}

\section{Summary}
\label{sec4}
We have developed a framework that can be used to study the
correlation properties of groups of Hilbert space indexed stochastic
processes. Our applications have been restricted to groups of size two
or three; however, it is clear that similar analyses are possible with
any finite number of processes. For example, the partial canonical
correlation work of Section~\ref{sec3} extends in principle to
examination of pairs of residual processes after correcting for
projections onto several other processes.

We note in passing that it has been assumed that all the
$\mathcal{H}$-valued random variables take values in the same Hilbert
space. The extension to where some or all of the variables produce
elements of different Hilbert spaces incurs some additional notational
expense but is otherwise straightforward.
\begin{appendix}
\section*{Technical Appendix}
\label{sec5}
In this Appendix, we collect some of the mathematical details that were
needed for our main results. In particular, the developments in
Sunder \cite{sunder} play a pivotal role in
Sections~\ref{sec2}--\ref{sec3}. Thus, we first summarize the key aspects
of that work that were employed in the paper.

Assume that a Hilbert space $\mathcal{H}$ can be written as the
algebraic direct sum of $n$ closed subspaces $M_1,\ldots,M_n$. That
is,
\begin{eqnarray*}
\mathcal{H}=\sum_{i=1}^nM_i,
\end{eqnarray*}
where $M_i\cap\sum_{j\neq i}M_j= \{0 \}$. Now, for $1\leq k
\leq n$ define
\begin{eqnarray*}
L_k = \Biggl( \sum_{i=1}^kM_i
\Biggr) \cap \Biggl( \sum_{i=1}^{k -
1}M_i
\Biggr)^\bot.
\end{eqnarray*}
Then, $L_k \bot M_i$, for $i=1,\ldots,k-1$, and by construction
$\sum_{i=1}^{k}L_i=\sum_{i=1}^kM_i$ for $k = 1, \ldots,n$.

Let $P_{M_k}$ and $P_{L_k}$ be the orthogonal projection operators
onto $M_k$ and $L_k$, respectively. Then, for $1\leq k\leq n$ and
$1\leq j \leq k\leq n$ we define the restriction of $P_{L_j}$ to $M_k$
by $P_{L_j|M_k}x=P_{L_j}x$ for $x \in M_k$ and use
$P_{M_k|L_j}y=P_{M_k}y$ for $y \in L_j$ to indicate the restriction of
$P_{M_k}$ to $L_j$. Sunder \cite{sunder} establishes the following
relationship between the $M_k$ and $L_k$.

%
\begin{theorem}
\label{thsunder}
For $x\in M_k$, we can write $M_k$ as
\begin{eqnarray*}
M_k & =& \bigl\{(P_{L_1|M_k}x,\ldots,P_{L_{k-1}|M_k}x,P_{L_k|M_k}x,0,
\ldots ,0) \bigr\}
\\
&=& \bigl\{\bigl(P_{L_1|M_k}(P_{L_k|M_k})^{-1}P_{L_k|M_k}x,
\ldots, P_{L_k|M_k}x,0,\ldots,0\bigr) \bigr\}
\\
&=& \bigl\{(A_{L_1|L_k}z,\ldots,A_{L_{k-1}|L_k}z,z,0,\ldots,0) \bigr\},
\end{eqnarray*}
where $z=P_{L_k|M_k}x\in L_k$ and
$A_{L_j|L_k}=P_{L_j|M_k}(P_{L_k|M_k})^{-1}$ for $1\leq j\leq k \leq
n$.
\end{theorem}

Theorem~\ref{thsunder} has the consequence that problems involving
optimization over $M_k$ can instead be formulated in terms of
equivalent problems on $L_k$ which is how it is applied in
Sections~\ref{sec2}--\ref{sec3}.

We next turn to the proof of Theorem~\ref{thccachar}. This is
accomplished via the following proposition and its corollaries.

%
\begin{proposition}
\label{proPh}
If $h=(C_{12}f_2,f_2) \in M_2$, then
$P_{L_1|M_2}h=(C_{12}f_2,C_{21}C_{12}f_2)$ and
$P_{L_2|M_2}h=(I-P_{L_1|M_2})h=(0,C_{22.1}f_2)$.
\end{proposition}

\begin{pf} Let $h_1=(f_1,C_{21}f_1)\in M_1 =L_1$. Then,
\begin{eqnarray*}
\langle P_{L_1|M_2}h_2,h_1\rangle _{\mathcal{H}(Q)}=
\langle h_2,h_1\rangle _{\mathcal{H}(Q)}
\end{eqnarray*}
for every $h_1\in M_1$. Writing $P_{L_1|M_2}h_2=(f_1^\star
,C_{21}f_1^\star)$ leads to
\begin{eqnarray*}
\langle P_{L_1|M_2}h_2,h_1\rangle
_{\mathcal{H}(Q)}&=&\bigl\langle \bigl(f_1^\star,C_{21}f_1^\star
\bigr),(f_1,0)\bigr\rangle _0 = \bigl\langle
f_1^\star,f_1\bigr\rangle _{\mathcal{H}(S_1)}
\\
&=& \bigl\langle (C_{12}f_2,f_2),h_1
\bigr\rangle _{\mathcal{H}(Q)} = \bigl\langle (C_{12}f_2,f_2),(f_1,0)
\bigr\rangle _0
\\
&=& \langle C_{12}f_2,f_1\rangle
_{\mathcal{H}(S_1)}
\end{eqnarray*}
for every $f_1\in\mathcal{H}(S_1)$. So, $f_1^\star=C_{12}f_2$.
\end{pf}

\begin{corollary}
If $h=(0,\tilde{f}_2)\in L_2$,
$(P_{L_2|M_2})^{-1}h=(C_{12}C_{22.1}^{-1}\tilde
{f_2},C_{22.1}^{-1}\tilde{f_2})$.
\end{corollary}

%
\begin{corollary}
\label{corBh}
For $h=(0,\tilde{f_2})\in L_2$, we have
\begin{eqnarray*}
Bh:=P_{L_1|M_2}(P_{L_2|M_2})^{-1}h=\bigl(C_{12}C_{22.1}^{-1}
\tilde{f_2}, C_{21}C_{12}C_{22.1}^{-1}
\tilde{f_2}\bigr).
\end{eqnarray*}
\end{corollary}

%
\begin{corollary}
\label{corh}
Let $h = (0,\tilde{f_2}), h' = (0,\tilde{f}_2') \in L_2$. Then,
\begin{eqnarray*}
\bigl\langle h,h'\bigr\rangle _{\mathcal{H}(Q)}=\bigl\langle (0,
\tilde{f_2}),Q^{-1}\bigl(0,\tilde {f}_2'
\bigr)\bigr\rangle _0=\bigl\langle \tilde{f_2},
C_{22.1}^{-1}\tilde{f}_2'\bigr\rangle
_{\mathcal{H}(S_2)}.
\end{eqnarray*}
\end{corollary}

With a little extra effort we also obtain the following corollary.

%
\begin{corollary}
\label{corB*}
$B^*(f_1,C_{21}f_1) =(0,C_{21}f_1)$.
\end{corollary}

\begin{pf} For $h=(f_2,C_{21}f_1)\in M_1=L_1$ and
$\tilde{h}=(0,\tilde{f_2})\in L_2$,
\begin{eqnarray*}
\langle h,B\tilde{h}\rangle _{\mathcal{H}(Q)}&=&\bigl\langle Q^{-1}h,B
\tilde{h}\bigr\rangle _0
\\
&=& \bigl\langle (f_1,0),\bigl(C_{12}C_{22.1}^{-1}
\tilde {f_2},C_{21}C_{12}C_{22.1}^{-1}
\tilde {f_2}\bigr)\bigr\rangle _0
\\
&=& \bigl\langle f_1,C_{12}C_{22.1}^{-1}
\tilde{f_2}\bigr\rangle _{\mathcal{H}(S_1)} = \bigl\langle
C_{22.1}^{-1}C_{21}f_1,
\tilde{f_2}\bigr\rangle _{\mathcal{H}(S_2)}
\\
&=& \bigl\langle B^*h,\tilde{h}\bigr\rangle _{\mathcal{H}(Q)}.
\end{eqnarray*}
An application of Corollary~\ref{corh} completes the proof.
\end{pf}

Finally, we give the details for proving Theorem~\ref{thmpccabb}. Analogous to the proof of Theorem~\ref{thccachar},
the steps are broken down into a proposition and its subsequent
corollaries.

%
\begin{proposition}
\label{proPPh}
If $h=(C_{12}f_2,f_2,C_{32}f_2)$,
$P_{L_1|M_2}h=(C_{12}f_2,C_{21}C_{12}f_2,C_{31}C_{12}f_2)$
and
$P_{L_2|M_2}h=(I-P_{L_1|M_2})h=(0,
C_{22.1}f_2,(C_{32}-C_{31}C_{12})f_2)$.
\end{proposition}

\begin{pf} For $h_1=(f_1,C_{21}f_1,C_{31}f_1) \in M_1=L_1$,
we have the relation
\begin{eqnarray*}
\langle P_{L_1|M_2}h,h_1 \rangle_{\mathcal{H}(Q)} =\langle
h,h_1 \rangle_{\mathcal{H}(Q)}.
\end{eqnarray*}
Writing $P_{L_1|M_2}h=(f_1^\star,C_{21}f_1^\star,C_{31}f_1^\star)$
leads to
\begin{eqnarray*}
\langle P_{L_1|M_2}h,h_1 \rangle_{\mathcal{H}(Q)}&=&\bigl\langle
\bigl(f_1^\star,C_{21}f_1^\star,C_{31}f_1^\star
\bigr),(f_1,0,0) \bigr\rangle_0
\\
&=& \bigl\langle f_1^\star, f_1 \bigr
\rangle_{\mathcal{H}(S_1)}
\\
&=& \langle h, h_1 \rangle_{\mathcal{H}(Q)}
\\
&=& \bigl\langle(C_{12}f_2,f_2,C_{32}f_2),
(f_1,0,0) \bigr\rangle_0
\\
&=& \langle C_{12}f_2,f_1
\rangle_{\mathcal{H}(S_1)}
\end{eqnarray*}
for every $f_i \in\mathcal{H}(S_i)$ with $i=1,2$. So
$f_1^\star=C_{12}f_2$.
\end{pf}

For subsequent notational convenience, let
\begin{eqnarray*}
C_0=C_{33.1}-(C_{32}-C_{31}C_{12})C_{22.1}^{-1}(C_{23}-C_{21}C_{13}).
\end{eqnarray*}

%
\begin{corollary}
\label{corl3form}
If $h=(C_{13}f_3,C_{23}f_3,f_3)$,
$P_{L_1|M_3}h=(C_{13}f_3,C_{21}C_{13}f_3,C_{31}C_{13}f_3)$,\linebreak[4]
$P_{L_2|M_3}h=(0,(C_{23}-C_{21}C_{13})f_3,(C_{33.1}-C_0)f_3)$ and
$P_{L_3|M_3}h=(0,0,C_0f_3)$.
\end{corollary}

\begin{pf} For $\tilde{h}_2=(0,
C_{22.1}f_2,(C_{32}-C_{31}C_{12})f_2) \in L_2$ and $h \in M_3$, we
have the\vspace*{1pt} relation $\langle P_{L_2|M_3}h, \tilde{h}_2
\rangle_{\mathcal{H}(Q)}=\langle h,\tilde{h}_2
\rangle_{\mathcal{H}(Q)}$. If we write $P_{L_2|M_3}h=(0,
C_{22.1}f_2^\star,(C_{32}-C_{31}C{12})f_2^\star)$, then
\begin{eqnarray*}
&&\langle P_{L_2|M_3}h, \tilde{h}_2 \rangle_{\mathcal{H}(Q)}
\\
&&\quad = \bigl\langle\bigl(0, C_{22.1}f_2^\star,(C_{32}-C_{31}C{12})f_2^\star
\bigr), (-C_{12}f_2,f_2,0) \bigr
\rangle_0
\\
&&\quad  = \bigl\langle C_{22.1}f_2^\star,
f_2 \bigr\rangle_{\mathcal{H}(S_2)}
\\
&& \quad = \langle h,\tilde{h}_2 \rangle_{\mathcal{H}(Q)}
\\
&&\quad  = \bigl\langle(0,0,f_3), \bigl(0,C_{22.1}f_2,(C_{32}-C_{31}C_{12})f_2
\bigr) \bigr\rangle_0
\\
&&\quad  = \bigl\langle f_3, (C_{32}-C_{31}C_{12})f_2
\bigr\rangle_{\mathcal{H}(S_3)}
\\
&&\quad  = \bigl\langle (C_{23}-C_{21}C_{13})f_3,
f_2 \bigr\rangle_{\mathcal{H}(S_2)}.
\end{eqnarray*}
So, $f_2^\star=C_{22.1}^{-1}(C_{23}-C_{21}C_{13})f_3$.
\end{pf}

\begin{corollary}
\label{proPBh}
For $h=(0,0,\tilde{f}_3) \in L_3$,
\begin{eqnarray*}
Bh = \bigl(0,(C_{23}-C_{21}C_{13})C_0^{-1}
\tilde{f}_3, (C_{32}-C_{31}C_{12})C_{22.1}^{-1}(C_{23}-C_{21}C_{13})C_0^{-1}
\tilde{f}_3\bigr).
\end{eqnarray*}
\end{corollary}

%
\begin{corollary}
\label{proPB*h}
If $h=(0, C_{22.1}f_2,(C_{32}-C_{31}C_{12})f_2) \in L_2$, then
\begin{eqnarray*}
B^*h=\bigl(0,0,(C_{32}-C_{31}C_{12})f_2
\bigr).
\end{eqnarray*}
\end{corollary}

\begin{pf} For $h=(0, C_{22.1}f_2,(C_{32}-C_{31}C_{12})f_2) \in
L_2$ and $\tilde{h}_3=(0,0,\tilde{f}_3) \in L_3$,
\begin{eqnarray*}
\langle B\tilde{h}_3, h \rangle_{\mathcal{H}(Q)} &=& \bigl\langle B
\tilde{h}_3, Q^{-1}h \bigr\rangle_0
\\
&=& \bigl\langle B\tilde{h}_3, (-C_{12}f_2,f_2,0)
\bigr\rangle_0
\\
&=& \bigl\langle (C_{23}-C_{21}C_{13})C_0^{-1}
\tilde{f}_3, f_2 \bigr\rangle_{\mathcal{H}(S_2)}
\\
&=& \bigl\langle C_0^{-1}\tilde{f}_3,
(C_{32}-C_{31}C_{12})f_2 \bigr
\rangle_{\mathcal{H}(S_3)}
\\
&=& \bigl\langle \tilde{h}_3, B^*h \bigr\rangle_{\mathcal{H}(Q)}
\\
&=& \bigl\langle Q^{-1}\tilde{h}_3,B^*h \bigr
\rangle_0
\\
&=& \bigl\langle \bigl(\bigl[C_{21}C_{22.1}^{-1}(C_{23}-C_{21}C_{13})-C_{13}
\bigr]C_0^{-1}\tilde {f}_3,
\\
&&\hphantom{ \bigl\langle \bigl(}{} -C_{22.1}^{-1}(C_{23}-C_{21}C_{13})C_0^{-1}
\tilde{f}_3, C_0^{-1}\tilde{f}_3
\bigr), B^*h \bigr\rangle_0.
\end{eqnarray*}
\upqed
\end{pf}
\end{appendix}

\section*{Acknowledgements}
The authors' research was supported by NSF Grant DMS 0652833.
Rosemary Renaut acknowledges the support of AFOSR Grant 025717:
Development and Analysis of Non-Classical Numerical Approximation
Methods, and NSF Grant DMS 1216559: Novel Numerical Approximation
Techniques for Non-Standard Sampling Regimes.
A helpful consultation
with Randy Eubank and Jack Spielberg, and the input of two referees and
an associate
editor are gratefully acknowledged.




\printhistory

\end{document}